\documentclass[onecolumn,notitlepage,aps,pra,10pt]{revtex4-1}
\usepackage{amsmath}
\usepackage{color}
\usepackage{graphicx} 

\begin{document}

\title{Tight Bounds for Tight Links: Ropelength of $T(Q,Q)$ torus links}
\author{Alexander R. Klotz}
\affiliation{Department of Physics and Astronomy, California State University, Long Beach}

\begin{abstract}
    Ropelength, $L$, is a parameter characterizing the minimum contour length of a knot or link. There exist upper and lower bounds on ropelength with respect to crossing number, $C$, including a universal lower bound constraining $L\geq\alpha_0 C^{3/4}$ for some constant $\alpha_0$. There is currently an order-of-magnitude range for the value of $\alpha_0$ between 1.105 and 10.76. In this work, we show that $T(Q,Q)$ torus links can be constructed such that the upper bound is within a factor of 1.77 of the lower bound. We derive a stronger lower bound based on the convex hull around close-packed disks of approximately $\alpha_{T_{QQ}}>\sqrt{8\pi\sqrt{3}}+(2\pi+\sqrt{2\pi+7\sqrt{3}-12}\ )Q^{-1/2}\approx6.60+7.61Q^{-1/2}$, significantly higher than the best universal lower bound of 1.105. We show that a link can be constructed without any free parameters or geometric optimization that, when $Q$ is large, has a coefficient $\alpha_{T_{QQ}}<1.005\cdot 4\pi(5\sqrt{5}-8)/3\approx13.39$, and can be improved to to 11.68 by solving a helical no-overlap constraint equation that requires a conjectural approximation. For $Q$ up to 20 we construct links from smooth planar curves or toroidal helices minimized with respect to a small number of geometric parameters, that are between 6 and 60\% greater in ropelength than the lower bound. Many such links can be annealed to within 10\% of the lower bound using gradient descent. This represents significant progress towards developing sharp bounds on the ropelengths of specific classes of knots and links.
\end{abstract}

\maketitle

\section{Introduction}
Ropelength is a property characterizing the tightest possible embedding of a knot or link, typically expressed as the minimum contour length of a knot surrounded by a non-overlapping unit-thickness tube \cite{cantarella2002minimum}. Ropelength can be related to other knot invariants such as the braid index \cite{diao2020braid}, is known to be correlated with the determinant of a knot \cite{klotz2024ropelength}, and generally increases with respect to minimum crossing number. The tightest form of a given knot, known as an ideal knot, has a quasi-quantized space writhe \cite{pieranski2001quasi} that has strong but not fully understood relationship with the invariant signature of the knot \cite{klotz2024ropelength}. Ideal knots represent a useful system for studying the relationship between geometric and topological properties of knots. The author took an interest in ropelength to better understand the topology of knots that were arising stochastically in experiments with DNA molecules, which could only be distinguished by their apparent size. The exact ropelength is only known for a few constructable links, such as chains of Hopf links. There exist upper \cite{diao2019ropelengths} and lower \cite{diao2003lower} bounds for the ropelength of a knot with a given crossing number, and the value can be estimated for individual knots through numerical tightening algorithms \cite{ashton2011knot, pieranski1998search}. Typically these upper and lower bounds are quite far apart and not constraining given numerical estimates \cite{klotz2021ropelength}, but recent effort has been made to bring them closer together. 

It was recently shown that alternating knots have a ropelength $L$ that is at least linear with respect to crossing number $C$, with a lower bound of at least $L>C/56$ \cite{diao2022ropelength}. The tightest known alternating knots are $T(P,2)$ torus knots and links, which can be constructed from repeating helical units and are known to have a ropelength of at most $L<7.31C$ \cite{kim2024efficiency}. Non-alternating knots, however, are typically tighter and are known to have a ropelength greater than a constant $\alpha_0$ times the three-quarter power of crossing number, $L>\alpha_0 C^{3/4}$, although a sharp value of $\alpha_0$ is unknown. The lower bound on $\alpha_0$ is $(11/4\pi)^{-3/4}\approx 1.105$ \cite{buck1998four}, although the author has conjectured it may be at least 3.22 \cite{klotzmini}. While the upper bound for the ropelength of all knots is a multiple $C\log^5(C)$ \cite{diao2019ropelengths}, non-alternating $T(P,Q)$ torus knots have been shown to have an upper bound on ropelength that is also a three-quarter power of crossing number \cite{cantarella1998tight, diao1998complexity}, meaning the ropelength can be constrained by bringing the upper and lower three-quarter coefficients closer to each other.

The overall goal of this work is to construct torus links for which the upper and lower ropelength bounds are as close as possible. The first part focuses on establishing a strong lower bound based on minimal convex hulls around unit disks. The second part focuses on developing a method to efficiently construct links that, as $Q$ tends to infinity, have an asymptotic upper bound that is below twice the lower bound. The third part focuses on the specific cases of $Q=3-20$ for which the lower bound is highest, examining how close we can get with a few construction methods. Throughout this work we may use the term ropelength to refer to the smallest length achieved by a specific construction, rather than the ultimate tightest value which in most cases is unknown.

Empirically, torus knots and links have the lowest ropelength at a given crossing number \cite{ashton2011knot}. We focus on links of the form $T(pQ,Q)$ where lowercase $p$ is a small natural number such that $P=pQ$ and $C=pQ(Q-1)$. When constructed, each component wraps around the central hole of the torus once before closing with itself. This allows us to avoid challenges with defining a non-interfering closure scheme for a torus knot, and lets us define a stronger lower bound as well as examine how ropelength depends on $p$ for large $Q$. In a previous work, we established principles to efficiently construct torus links from concentric helices wrapped into a torus \cite{thompson2025ropelength}. We primarily focused on efficiently constructing helices such that the ropelength per helical twist was minimized, which would allow the most efficient \textit{linear} ratio of $L/C$ when $p$ is large. When minimizing the ratio of $\alpha=L/C^{3/4}$ we showed that our construction was optimized when $p=3$ and showed that toroidal links with ropelengths of $L\approx 12C^{3/4}$ could be constructed up to $Q=39$, and that large links could be constructed with  $\alpha_{T_{3Q,Q}}<19.11$. We argued that the asymptotic lower bound for such a construction was approximately $\alpha_{T_{3Q,Q}}>5$, bringing the upper and lower bound within a factor of 4.

In this work, we expand our lower bound to all $p$ and $Q$ using the Wegner inequality \cite{wegner1986uber} for the convex hull of unit disks. We show that when $p=1$ and $Q$ is small, the lower bound is high enough that it is comparable to tight links that can be constructed, and is asymptotically $\sqrt{8\pi\sqrt{3}}\approx6.60$ when $Q$ is large. Focusing directly on toroidal helices instead of linear helices, we outline a procedure for constructing $T(Q,Q)$ links that have a limiting coefficient of 13.38 without requiring numerical optimization, and can reduce this to 11.68 by solving an approximation of the no-overlap constraint for cylindrical helices. Finally, we focus on constructing tight links by minimizing a small number of geometric parameters, for low $Q$ in the regime where the lower bound is high and compare these to the results of numerical tightening. We show for $Q=3,4,5$ that links constructed from gibbous planar curves are more efficient than toroidal helices, typically within 20\% of their lower bound. 

\section{Lower Bound}
\begin{figure}
    \centering
    \includegraphics[width=1\linewidth]{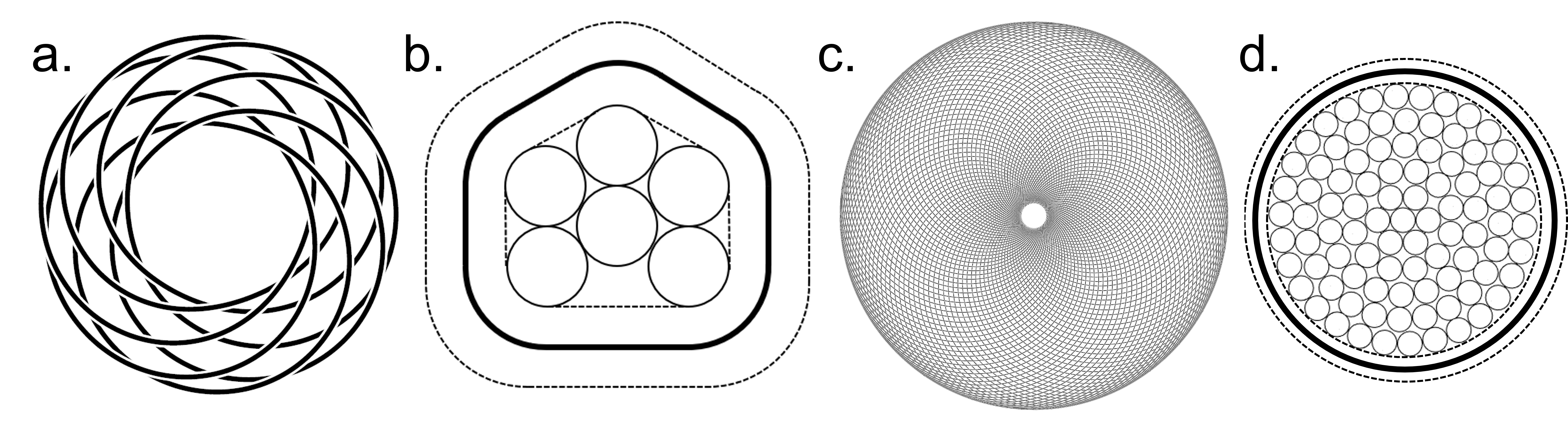}
    \caption{a. Diagram of a $T(7,7)$ link. b. The shape of the shortest curve that can be part of a unit-thickness $T(7,7)$ link, with its bold centerline pushed a distance of 1 from the minimal convex hull around 6 unit disks.   The lower bound arises from assuming that each component of a $T(Q,Q)$ link is the minimal convex hull of $Q-1$ unit disks. c. Diagram of a $T(92,92)$ link. d. An estimate of the minimal hull around 91 disks. For large enough $Q$ the minimal hull approximates a circular curve around a hexagonal packing.}
    \label{fig:lb}
\end{figure}

In a $T(pQ,Q)$ link, each component wraps $p$ times around the other $Q-1$ components. Our lower bound is based on treating each component in a $T(pQ,Q)$ link as the minimal convex hull around $p(Q-1)$ unit disks, as visualized in Fig. \ref{fig:lb}. The current universal lower bound for ropelength is based on constraining it to an electromagnetic knot energy \cite{buck1999thickness}, while the author's conjectured value constrains it to the Mobius energy \cite{klotzmini}. A stronger lower bound may be achieved for $T(pQ,Q)$ torus links. Following Cantarella et al. \cite{cantarella2002minimum}, their Theorem 10 states that a component linked with $N$ others must have a length greater than $2\pi$ plus the minimum convex hull around $N$ unit disks, where the extra $2\pi$ comes from pushing out from the convex hull by a unit distance. For $N=2,3,4,5$ the pushed minimal hulls have lengths of $4\pi+2N$, but for larger $N$ finding the length becomes non-trivial. In a $T(Q,Q)$ torus link, each of the $Q$ components is linked to the other $N=Q-1$, so that the ropelength must be at least:
\begin{equation}
    L_{T_{QQ}}>Q(4\pi+2(Q-1)), \ \ \ \ 3\leq Q\leq6
\end{equation}
This minimum length cannot be achieved without violating the no-overlap constraint, and beyond $N=5$ the minimal hull cannot link with $N$ copies of itself (e.g. Fig. 1b for $N=6$). More generally, if a component $J$ has a total linking number $Lk$ with other components in a link, its length must exceed the minimal hull around $Lk$ unit disks, plus $2\pi$. The isoperimetric inequality states that the minimal convex hull around a set of unit disks must have a greater perimeter than a circle with the same area as those disks, leading to Theorem 11 of Cantarella et al.:
\begin{equation}
    L_J>2\pi+2\pi\sqrt{Lk}
\end{equation}
In a $T(pQ,Q)$ link, each component has linking number $p$ with the others, giving it total $Lk=p(Q-1)$. This leads to the isoperimetric bound on the ropelength:
\begin{equation}
    L_{T_{pQ,Q}}>2\pi Q\left(1+\sqrt{p(Q-1)}\right)\approx2\pi\sqrt{p}Q^{3/2}
\end{equation}
Throughout this section, we will provide approximations for large $Q$. We are interested in the ratio  $\alpha=L/C^{3/4}$ of a given link, and use $\alpha_0$ to represent the universal lower bound and other subscripts for specific cases. Since the crossing number is $pQ(Q-1)$, the isoperimetric constraint, $\alpha_{iso}$ is:
\begin{equation}
    \alpha_{iso}=\frac{2\pi Q\left(1+\sqrt{p(Q-1)}\right)}{(pQ(Q-1))^{3/4}}\approx\frac{2\pi}{p^{1/4}}+\frac{2\pi}{p^{3/4}\sqrt{Q}}
\end{equation}
This is greatest for $p=1$ and when $Q$ is small, which is one reason we focus on $T(Q,Q)$ links. For $p=1$ and large $Q$, the limiting ratio is $2\pi$ which is a factor of 5.7 above the lower bound on $\alpha_0$. Because circles do not tile a plane, an isoperimetric hull around disks cannot be achieved, and when $Q$ is large, the minimal hull can be approximated a circle surrounding a hexagonal packing of disks with area fraction $\sigma=\pi/\sqrt{12}\approx0.91$. This has the effect of increasing the limiting value of $\alpha$ by about $\sigma^{-1/2}\approx 1.05$. The \textit{densest} hull around $N$ unit disks, known as a Wegner or Groemer packing \cite{boroczky2007note}, is that which best approximates a hexagonal shape surrounding a hexagonal packing. The minimal length hull will have an \textit{area} that is greater than or equal to a Wegner packing and a \textit{perimeter} that is greater than or equal to that of a circle with the same area. This leads to a stronger lower bound on the hull of $N$ hexagonally packed disks, given by the Wegner inequality \cite{wegner1986uber}:
\begin{equation}
    W(N)\geq\sqrt{4\pi\left(\sqrt{12}(N-1)+(2-\sqrt{3})\lceil\sqrt{12N-3}-3\rceil+\pi\right)}
\end{equation}
By pushing this by $2\pi$, incorporating the linking number for $N$, and multiplying it by $Q$ we arrive at our lower bound:
\begin{equation}
    L_{T_{pQQ}}\geq Q\left(2\pi+\sqrt{4\pi\left(\sqrt{12}(p(Q-1)-1)+(2-\sqrt{3})\lceil\sqrt{12p(Q-1)-3}-3\rceil+\pi\right)}\right)
\end{equation}
When $p=1$ and $Q$ is large, the dominant term is $2\pi/ Q^{3/2}/\sqrt{\sigma}$ as expected. We can find the coefficient $\alpha_w$ by dividing by the $3/4$ power of the crossing number:
\begin{equation}
        \alpha_w=\frac{ Q\left(2\pi+\sqrt{4\pi\left(\sqrt{12}(p(Q-1)-1)+(2-\sqrt{3})\lceil\sqrt{12p(Q-1)-3}-3\rceil+\pi\right)}\right)}{(pQ(Q-1)^{3/4}}
\end{equation}
This can be approximated by ignoring the ceiling function and taking a series expansion about large $Q$:
\begin{equation}
    \alpha_w>\sqrt{8\pi}\left(\frac{3}{p}\right)^{1/4} + \left(2\pi+\sqrt{2\pi(7\sqrt{3}-12)}\right)\frac{1}{\sqrt{Q}p^{3/4}}\approx\frac{6.60}{p^{1/4}}+\frac{7.61}{\sqrt{Q}p^{3/4}}
\end{equation}
The limiting value falls below the universal lower bound of 1.105 for $p>1271$, and below the author's conjectured bounds of 2 and 3.22 when $p>118$ and $p>17$ respectively. When $p=1$ the limiting value is $\sqrt{8\pi\sqrt{3}}\approx6.60$, six times the lower bound of $\alpha_0$. The bound falls below 10 for $Q>6$ and below 8 for $Q>28$. No investigation has produced a knot or link with $\alpha<10$, so the tightest known links must be made 40\% tighter to show this bound is sharp, which it likely is not. In our previous work we discussed this bound in the the specific case of $p=3$, for which efficient toroidal helices can be constructed and the Wegner coefficient is close to 5.  It is important to note that this bound may not apply when $p>1$ as the formal derivation of the lower bound involves a certain number of punctures of the surface of a cone by the link \cite{cantarella2002minimum}, a constraint which is not guaranteed to be satisfied. The isoperimetric bound should be treated as more rigorous in these cases, although the true lower bound is likely greater. 

The Wegner formula likely underestimates the minimal hull perimeter for finite $N$, and indeed it cannot be achieved for one in twenty-four numbers, beginning with 121 \cite{boroczky2007note}. Numerical estimates of the minimal convex hulls around unit disks \cite{kallrath2021near} suggest values that are on the order of 5\% greater than the Wegner formula, increasing the lower bound accordingly. Because these values are conjectural, we do not use them to generate a lower bound, but compare some our results to them in the final section.

\section{Zero-parameter Upper Bounds}
We desire a way to construct an efficient $T(Q,Q)$ link that is guaranteed to be non-overlapping and does not require any geometric or combinatorial optimization. Here we describe methods based on pairs of concentric toroidal helices that are linked through each others' donut holes. The first method will produce a link with a limiting ropelength of $13.38C^{3/4}$, which is always less than a factor of $2.04$ above the lower bound. The second method uses an approximation to the helical no-overlap constraint equation to create links with a limiting ropelength of $11.68C^{3/4}$. Our numerical parameters are typically precisely defined, but only apply in the limit of very large $Q$.

\subsection{Toroidal Helices}

Our constructions are based on concentric shells of cylindrical helices, wrapped into a torus. Here we will describe the basic features of toroidal helices. In our previous work, we established efficient ways to organize concentric helices such that in the large-$Q$ limit, the length of a single helical unit was below $7.83Q^{3/2}$, which establishes a strong linear ropelength bound for large $p$. The most efficient linear helix may not be the best toroidal helix, as the increase in the major radius of the torus to prevent overlap must also be taken into account.

If $N$ helices lie evenly spaced around a cylinder of radius $r$ and height $H$, they can repeat $p$ times and wrap into a toroidal $N$-helix of major radius $R_0=pH/2\pi$. The equation of the \textit{i}th helix is:
\begin{equation}
    x(\theta)=\left(R_0+r\cos(p\theta)\right)\cos\left(\theta+\frac{2\pi}{N}i\right), \ \ \ \     y(\theta)=\left(R_0+r\cos(p\theta)\right)\sin\left(\theta+\frac{2\pi}{N}i\right),\ \ \ \
z(\theta)=r\sin(p\theta)
\end{equation}
Here we are primarily concerned with $p=1$. The total length of a helix is $L=\sqrt{H^2+(2\pi r)^2}$, the length of a toroidal helix is larger by a small fraction (discussed below) that decreases as $H$ exceeds $r$. Two helices on the same cylinder are non-overlapping if the minimum distance between them is greater than or equal to 2. This can be described by the transcendental constraint equation \cite{huh2016best}:
\begin{equation}
    d=\sqrt{2r^2\left(1-\cos\left(\theta_\Delta+\frac{2\pi}{N}\right)\right)+\left(\frac{H\theta_\Delta}{2\pi}\right)^2}\geq2
\end{equation}
$-\pi\leq \theta_\Delta\leq\pi$ is the phase difference between the central point on helix $i$ and any point on helix $i-1$. There is a small negative value of $\theta_\Delta$ that minimizes this, but there is not a closed form solution for it. The challenge is to find parsimonious values of $r$, $H$, and $N$ that satisfy the constraint while minimizing helical length. For a single $N$-helix, the solution is that $r=\sqrt{2}N/\pi$ and $H=2\pi r=\sqrt{8}N$, assuming $N$ is large. When multiple concentric shells of helices are present, at most one shell will have this ``ideal'' configuration, and the challenge becomes choosing the arrangement of helices that will minimize the total length. Each shell must have the same height, and the radial separation between successive shells must be at least 2. We refer the reader to our previous work for a discussion on the best way to optimize this for a vertical helix \cite{thompson2025ropelength}.

If the radius of a helical shell is such that any reduction would cause an overlap, it is said to be radially constrained. Likewise, a system that cannot be reduced in height without overlap is said to be vertically constrained. The ideal $N$-helix is simultaneously radially and vertically constrained. In all other cases, a shell of helices is bound by at most one constraint. These constraints are modified slightly when discussing toroidal helices. When an $N$-helix is wrapped into a torus, the portions of the helix on the inside of the donut hole are brought into closer proximity, which will lead to overlaps if the straight helix is already constrained. To resolve this, we can increase the major radius $R_0$ by the minor radius. A system of concentric toroidal helices with an outer radius $r_o$ has an effective constraining height $H_c=2\pi(R_0-r_o)<H$ corresponding to the circumference of the donut hole. We define the inner radius of the donut hole $h=R_0-r_o$, using the symbol $h$ for consistency with previous work \cite{huh2016best}.

If two copies are made of a $T(Q,Q)$ link, and one copy is rotated and translated and then Hopf linked through the donut hole of the other, the result is a $T(2Q, 2Q)$ link with a length of at least $2L$ but $4C$ crossings, potentially reducing $\alpha$ by a factor of $\sqrt{2}$. These links are not torus links for $p>1$, as components will have linking number $p$ with members of the same torus and linking number $1$ with members of the other torus. In the case of $p=1$, if such a link is constructed from a torus link and its mirror image, the result is not the same. For example, linking two Hopf links in such a way yields a $T(4,4)$ which is equivalent to $L12n2206$, but linking a Hopf link with its reflection yields the slightly different $L12n2209$. This donut linking operation requires that $R_o\geq2r_o+2$.

\subsection{Constant Increment Construction}

\begin{figure}
    \centering
    \includegraphics[width=1\linewidth]{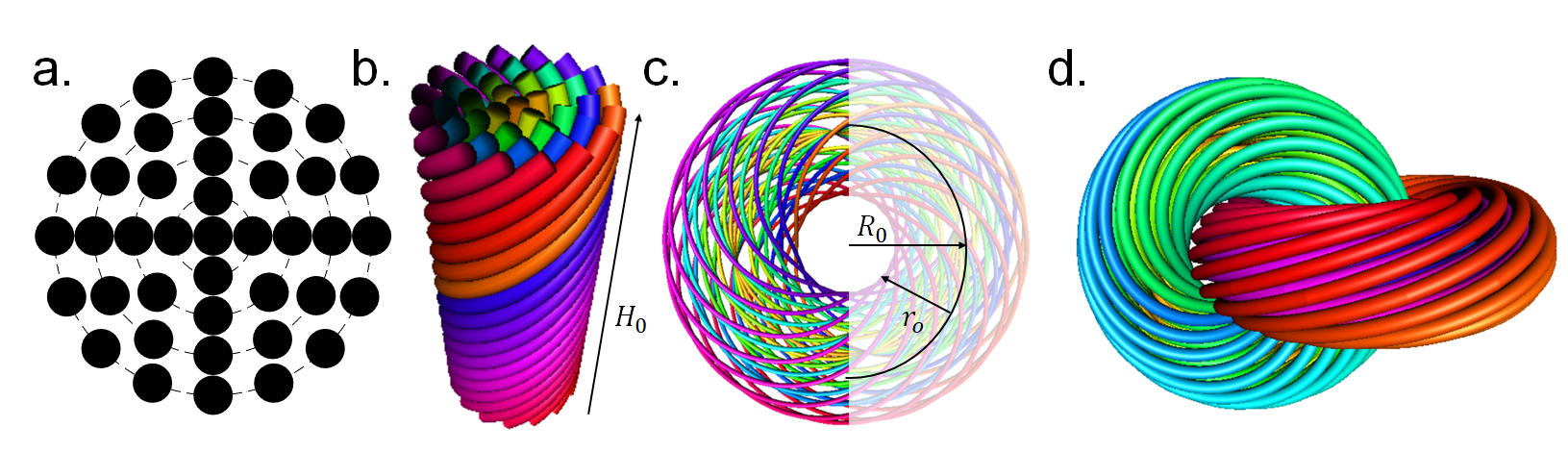}
    \caption{Construction of an efficient $T(Q,Q)$ link from concentric helices. a. The arrangement of helices, incrementing 4 per shell. b. The vertical concentric helices. c. A torus formed from wrapping the helix, shown with thinner curves. d. Two copies of the torus are increased in major radius and linked through each other, shown with unit thickness curves. The visible extra space highlights the unconstrained nature of the helices.}
    \label{fig:ub}
\end{figure}

We can construct a helical unit from a rod surrounded by concentric shells of $N$=4, 8, 12, 16, etc (Fig. \ref{fig:ub}a,b). helices at radii of  $r$=2, 4, 6, 8, etc. An $N$-helix has an ideal radius of approximately $0.45N$ such that $r_{ideal}\approx2.22N$. Every helix is radially underconstrained since each shell contains 0.9 times its ideal number of helices. With $T$ total shells, the outer shell has radius $r_o=2T$ and $N_o=4T$. The challenge is determining the minimum height of the helix such that when it is closed into a torus (Fig. \ref{fig:ub}c), the helices in the outer shell will be non-overlapping. In such a case, all the interior helices will be underconstrained in both radius and height and only the outer shell is constrained.

To determine that height we unwrap the cylinder containing the helices into a rectangle of height $H_c$ and width $2\pi r_o$, on which the helices form  lines parallel to the diagonal of the rectangle. Here were use the constrained height corresponding to the circumference of the donut hole, not the full height corresponding to the circumference of the centerline of the torus. Since the radius of each shell is half the number of helices it contains, we can write an expression for the pitch angle of each helix as $\tan\nu=\frac{H_c}{\pi N_o}$. The minimum distance between adjacent helices is their separation perpendicular to their pitch, $(H_c/N_o)\cos\nu$. We can constrain this minimum distance to be 2:

\begin{equation}
    \frac{H_c}{N_o}\cos(\arctan\frac{H_c}{\pi N_o})\geq2\rightarrow H_c \geq \frac{2\pi}{\sqrt{\pi^2-4}}N_o\approx2.59 N_o
    \label{eq:rect}
\end{equation}
This is approximately 92\% the ideal height of an $N$-helix, $\sqrt{8}N$. Our approximation using the perpendicular distance between two lines on a rectangle becomes less valid as the circumference of the cylinder becomes small compared to its height. In this case the circumference is slightly larger than the height. The validity of this approximation is discussed in the next subsection and the Appendix. Solving the transcendental constraint equation for the smallest possible radius yields a height requirement of $2.628N_o$, but these constructions are meant to apply to large systems. 

If the helices of height $H_c$ are closed into a torus with major radius $R=H_c/2\pi$, the extrema on the inside of the torus would overlap. To account for this, the major radius $R_0$ must be increased by the radius of the outer helices:
\begin{equation}
 R_0=\frac{H_c}{2\pi}+r_o=\frac{2\pi}{\sqrt{\pi^2-4}}\frac{N_o}{2\pi}+\frac{N_o}{2}=\left(\frac{1}{\sqrt{\pi^2-4}}+\frac{1}{2}\right)N_o\approx0.913N_o   
\end{equation}
This is sufficient to construct such a helix. To determine its ropelength in the large-$Q$ limit, we can set up a summation over the length of each shell of helices and evaluate it as an integral over all the shells. To do so it is useful to express all quantities in terms of $Q$. $Q$ can be found by summing multiples of 4 from from 4 to $4T$, giving $T\approx\sqrt{Q/2}$ and $N_0\approx\sqrt{8Q}$. 

The proper height of the helix that is closed to have that radius is:
\begin{equation}
H_0=2\pi R_0=\frac{H_c}{2\pi}+r_o=2\pi\left(\frac{1}{\sqrt{\pi^2-4}}+\frac{1}{2}\right)N_o=2\pi\left(\frac{1}{\sqrt{\pi^2-4}}+\frac{1}{2}\right)\sqrt{8Q}\approx16.22\sqrt{Q}\equiv\kappa\sqrt{Q}
\end{equation}
To compute the sum we assume $Q$ is large and ignore terms that vanish in that limit. The contour length of all the helices in the torus is the same as that of the helix, which we can write as:
\begin{equation}
    L_4=\sum_{i=1}^T4i\sqrt{H_0^2+(2\pi r_i)^2}=\sum_{i=1}^{\sqrt{Q/2}}4i\sqrt{\kappa^2 Q+(4\pi i)^2 }
\end{equation}
The integral approximation is:
\begin{equation}
    L_4\approx\int_{i=1}^{\sqrt{Q/2}}4i\sqrt{\kappa^2    Q+(4\pi i)^2 }di=\frac{Q^{3/2}}{12\pi^2}\left((\kappa^2+8\pi^2)^{3/2}-\left(\kappa^2+\frac{16\pi^2}{Q}\right)^{3/2} \right)
\end{equation}
When $Q$ is large, the term with $1/Q$ vanishes and the crossing number is approximately $Q^2$, giving $C^{3/4}\approx Q^{3/2}$ and a three-quarter coefficient of:
\begin{equation}
    \alpha_4=\frac{1}{12\pi^2}\left((\kappa^2+8\pi^2)^{3/2}-\kappa^3 \right)\approx17.36
\end{equation}
This is a factor of 2.63 above the lower bound and is an improvement on the construction described in our previous work \cite{thompson2025ropelength}, but can be significantly improved. We note that when $p=1$ any of the outer helices may be smoothed into a circle with a normal perpendicular to the innermost circular helix, reducing the ropelength sightly. This may be done with a second helix, which will not be circular but will still be reduced in length, and then a third, etc, until the link is essentially two separate $T(Q/2,Q/2)$ links Hopf linked with each other. Likewise, a $T(Q,Q)$ link can be doubled and linked through its own hole with a copy to form a $T(2Q,2Q)$ link (Fig. \ref{fig:ub}d). Applying this copying and linking procedure to a $T(pQ,Q)$ link will double the length and increase the crossing number to $2pQ(Q-1)+2Q^2$, the $2Q^2$ arising from the fact that each component of each torus has $Lk=1$ with each component in the other. When $p=1$ the crossing number is quadrupled and $\alpha$ can decrease only if the major radius does not increase by a factor of more than $\sqrt{2}$, but this procedure will not reduce $\alpha$ if $p>3$. When $p>1$, the new construction is no longer a $T(pQ,Q)$ torus link but is subject the lower bound of a $T((p+1)Q,Q)$ link. 

In order for one copy of the torus to fit through another, its donut hole must be greater than the diameter of outer layer of helices, plus 2. When $Q$ is large this requires the major radius of the 4-incremented torus to increase by a factor of $\approx1.095$ to $4T=\sqrt{8Q}$ and the height term in the integral becomes $2\pi\sqrt{8}>\kappa$. Re-computing the limiting ropelength ratio from the integral we have:
\begin{equation}
    \alpha_{4D}=\frac{4\pi}{3}(5\sqrt{5}-8)\approx13.32
\end{equation}
This would be a significant improvement and is at most a factor of 2.019 above the lower bound. However, a toroidal helix is slightly longer than the equivalent straight helix, by a ratio that decreases with the ratio of major to minor radii. This excess length can be computed as an integral with no closed-form solution. The major radius in our first construction is approximately 1.81 times the minor radius of the outermost helices, with an excess length of 1.5\%. When extended to form a double torus, the major radius is twice the minor and the increase of the largest helix is 1.11\%. The rest of the helices have a smaller and vanishing length correction. To compute the integrated effect of this toroidal lengthening, we numerically compute the following limit, which is the ratio of the total length of toroidal helices to equivalent straight helices, and converges quickly.
\begin{equation}
\lim_{T\rightarrow\infty} \frac{\sum_{i=1}^Ti\int_0^{2\pi}\sqrt{4i^2+(4 T-2i\cos\theta)^2}d\theta}{\sum_{i=1}^T4\pi i\sqrt{4T^2+i^2}}\approx1.0042
\end{equation}
With this correction, the ropelength ratio converges to 13.38. If a firm upper bound is desired, the correction falls below 1.005 when $T\geq10$. Other constructions will have other versions of this limit. For the single 4-incremented helix, the correction is 0.58\%. The helical corrections are the only numerical parameters in our derivations that we do not have exact expressions for; a table of such corrections appears in the Appendix.

Having a single torus with increments of 5 helices per shell instead of 4 is asymptotically less efficient (approaching a coefficient of $\alpha_5=19.2$). Because the tori can fit through their own donut holes when $T>4$, when linked together they approach a coefficient of 13.59, with a correction of 0.4\% to $\alpha_{5D}=13.64$. Both the 4- and 5-incremented constructions are less efficient when the outer shells are not full, as the outer helices increase the major radius of their tori. This can be resolved by only increasing the radius when the number of helices in the outer shell exceeds that of the penultimate shell, before which they are unconstrained. This correction improves the efficiency of the single 4-incremented torus, but the radius still needs to increased to create a double torus. Because the 5-incremented construction already has a major radius that is more than twice its minor radius, this correction significantly improves efficiency for $T>8$, and the 5-incremented helix is more efficient than the 4-incremented until $Q\approx23800$ or 566 million crossings. This is below twice the lower bound until $Q\approx15000$ and 225 million crossings, or $Q\approx2380$ if only the isoperimetric inequality is assumed. However, the vast majority of the helices in these systems are underconstrained and they can be made more efficient.

\subsection{Optimal solution}
The methods described above will produce an efficient torus link without requiring numerical optimization or constraint-checking. More efficient large-$Q$ helices can be created by placing the maximum number of helices that can fit in each shell without violating the no-overlap constraint. Removing helices from the outer shells and placing them in the inner shells, what we previously termed ``reverse Jenga,'' may reduce the length of a given construction. The 4-incremented torus, for example, requires only about 8.7\% of its shells or 16.7\% of its helices to be moved to an inner layer in order to fit through its own donut hole without increasing in length. Achieving this can reduce the limiting ratio to below 12.27, with an additional improvement due to the reverse-Jenga reduction. This is feasible as the lower 91\% of shells can accommodate an additional 21\% more helices.

We believe the ideal configuration has three features: each shell has the smallest possible minor radius, each torus has the smallest possible major radius that can fit through its own donut hole, and each shell has the most possible helices. The first condition is satisfied by having the \textit{i}th shell at radius $2i$, the second by having the major radius be twice the minor radius of the $T$th shell (that is, $4T$), but the third constraint is more challenging.

The number of helices that can fit in a given shell requires minimizing the transcendental constraint equation and solving for $N$ such that the minimum is not less than 2. We cannot do this symbolically, but a small angle approximation exists when $N$ is large. This is equivalent to finding the distance between two parallel lines on a rectangle, which we used in Eq. \ref{eq:rect}. This allows an approximate expression for the maximum number of helices in a shell:

\begin{equation}
    N_{a}=\frac{h\pi r}{\sqrt{h^2+r^2}}
\end{equation}
Where $h$ is the radius of the donut hole, taking the value of $2T$ in this case. This approximation is a slight overestimate for the number of helices that can fit in a shell, with the discrepancy decreasing as the radius of the shell approaches the effective height and both become large. In the Appendix we sketch a partial proof, based on a cubic approximation of the constraint equation, that there is a value $\epsilon$ that can be subtracted from $N_a$ such that the constraint is satisfied:
\begin{equation}
    N_{max}=\lfloor N_a-\epsilon\rfloor,\ \ \ \ \ \ \ \ \ \ \ \ \ \pi\left(r-\left(\arcsin\frac{1}{r}\right)^{-1}\right)\leq\epsilon\leq2\pi-6
\end{equation}
The correction is largest for an asymptotically tall helix of radius 2, where $N_a=2\pi$ but the maximum is 6. It is vanishing for larger $r$ particularly when rounding down to the nearest integer. The subtraction required from the approximate number per shell does not affect the asymptotic $\alpha$ of this construction. To compute the limiting ratio, we again set up a sum to compute the length and evaluate it as an integral. In this case we choose the maximum number of shells and additionally compute $Q$. Note that the number of helices per shell is computer from the constraint height, but the length is computed from the absolute height.

\begin{equation}
    L_{opt}=\sum_{i=1}^{T}\frac{4\pi iT}{\sqrt{4T^2+4i^2}}\sqrt{(8\pi T)^2+(4\pi i)^2}\approx\int_{1}^{T}8\pi^2iT\sqrt{\frac{4T^2+i^2}{T^2+i^2}}di=2 \pi ^2 \left[\sqrt{40}-\log (27)+6 \tanh ^{-1}\left(\sqrt{\frac{2}{5}}\right)-4\right]T^3
\end{equation}
This evaluation assumes $T$ is large, and the cubic coefficient is approximately 69. The number of helices evaluates to:
\begin{equation}
    Q=\sum_{i=1}^{T}\frac{2\pi iT}{\sqrt{i^2+T^2}}\approx2\pi T\left(\sqrt{2}T-\sqrt{1+T^2}\right)\approx2\pi(\sqrt{2}-1) T^2\approx 2.6T^2
\end{equation}
The coefficient of 2.6 implies a slightly higher density of helices than the 2.5 for the 5-increment system, which also requires a greater major radius. Taking the limit of $L/C^{3/4}$ as $Q$ tends to infinity and dividing by $\sqrt{2}$ to account for the donut linking procedure yields the ratio:

\begin{equation}
    \alpha_{opt}=\frac{L_{opt}}{C^{3/4}}=\sqrt{(7+5\sqrt{2})\pi}\left(\sqrt{10}-2+3\left(\tanh ^{-1}\sqrt{\frac{2}{5}}-\tanh ^{-1}\frac{1}{2}\right)\right)\approx11.64
\end{equation}
The toroidal correction in this case is 1.00385, briging the upper bound for large links to 11.68. Since the outer shell has $R_0=2r_o$, the correction 1.011 may be used as a universal overestimate. The upper bound is 1.77 times the Wegner bound and 1.86 times the isoperimetric bound. We believe this is the best coefficient that can be achieved with a concentric helical construction, and requires the assumption that the cubic approximation to the transcendental constraint equation does not yield overlapping configurations. This assumption is supported by numerical tests up to $T=100$ ($Q=52203$, $C\approx2.7\times10^9$) tested with 1000 points per helix. While this construction only covers specific values of $Q\approx5.2T^2$, outer helices can be removed to reach a desired $Q$.

Figure \ref{fig:ub} shows the coefficients from the doubled 4-increment, 5-increment, and optimal constructions as a function of $Q$. These lengths were computed by constructing helices according to each specification, and tested for overlaps when computationally feasible. The matching lines in red and black show the best- and worst-case scenarios of the incremented models, the best case scenario being that $4T-4$ or $5T-5$ helices were added to the outer row so that the major radius need not be increased, the worst case scenario being a full outer shell. All curves approach their limiting values as $Q$ becomes large, shown without the helical correction. There are noticeable discontinuities in some of the curves, which correspond to critical values of $T$ below which the major radius must be increased to allow for double linking. For the 4-incremented helix the positive trend in the worst case scenario below $Q=100$ corresponds to a single torus being more efficient than two. The same data is presented as a ratio with respect to the Wegner lower bound. The links described in this section can also be constructed if one of the copies of the torus is reflected, although the result is no longer a torus link.

 

\begin{figure}
    \centering
    \includegraphics[width=1\linewidth]{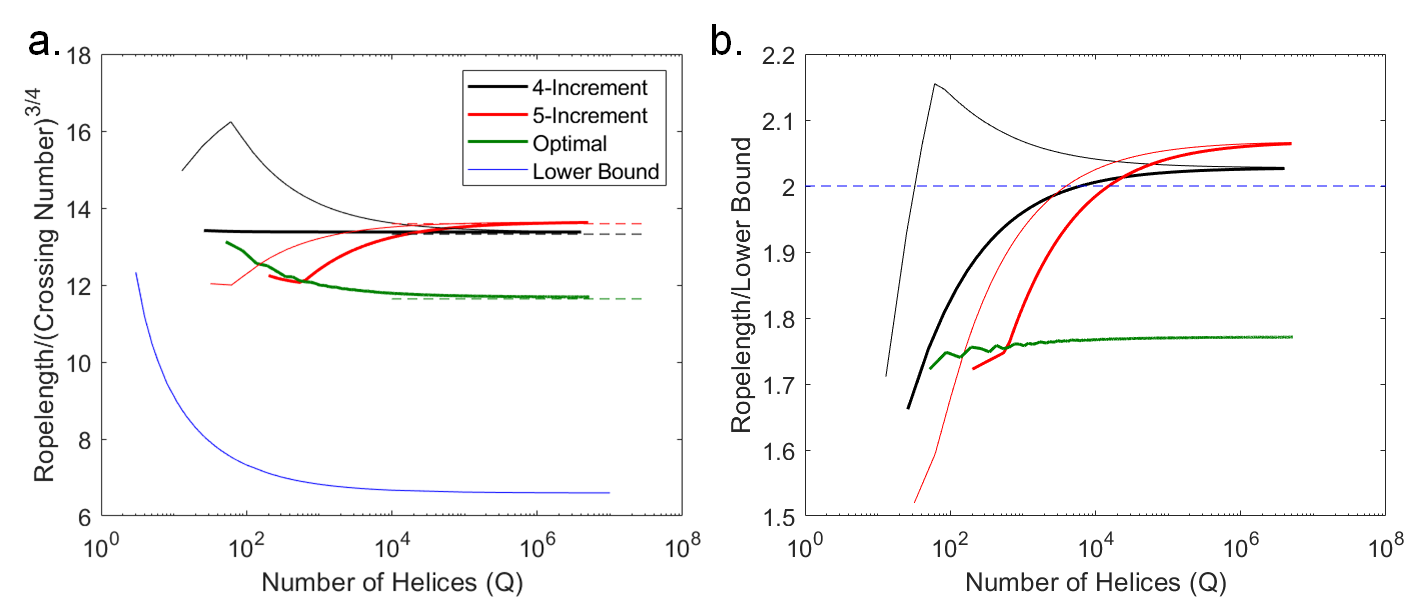}
    \caption{a. Three-quarter coefficients for the $T(Q,Q)$ links formed by doubly-linked tori of concentric helices incremented by 4 and 5 per shell, as well as the optimal construction. The two curves in each color show the best- and worst-case scenarios, which depend on how filled the outer helices are. Kinks correspond to configurations below which the major radius must be increased to allow for double torus linking. Limiting values are shown as dashed lines without the helical correction and the lower bound is shown for comparison. b. Shows the same data as a ratio relative to the lower bound.}
    \label{fig:ubgraph}
\end{figure}




\section{Low-$Q$ Constructions and Optimization}
Although we have shown that we can construct a helix that is within a factor of 1.77 of the lower bound when the crossing number is large, the lower bounds are stronger when $Q$ is small. Combinatorial explosion makes developing an understanding of knot invariants difficult as the crossing number increases, but there is only one $T(Q,Q)$ link for each $Q$, and it is worth putting effort into a finite number of links for which a configuration may be constructed that is close to the lower bound. Our initial goal was to determine whether geometries besides toroidal helices could be used to construct links close to the lower bounds using a minimal amount of numerical optimization. We were also interested in determining whether our best-constructed links have a ropelength that continues to decrease with increasing $Q$ as the lower bound does, or reaches a plateau independent of $Q$.

\begin{figure}
    \centering
    \includegraphics[width=1\linewidth]{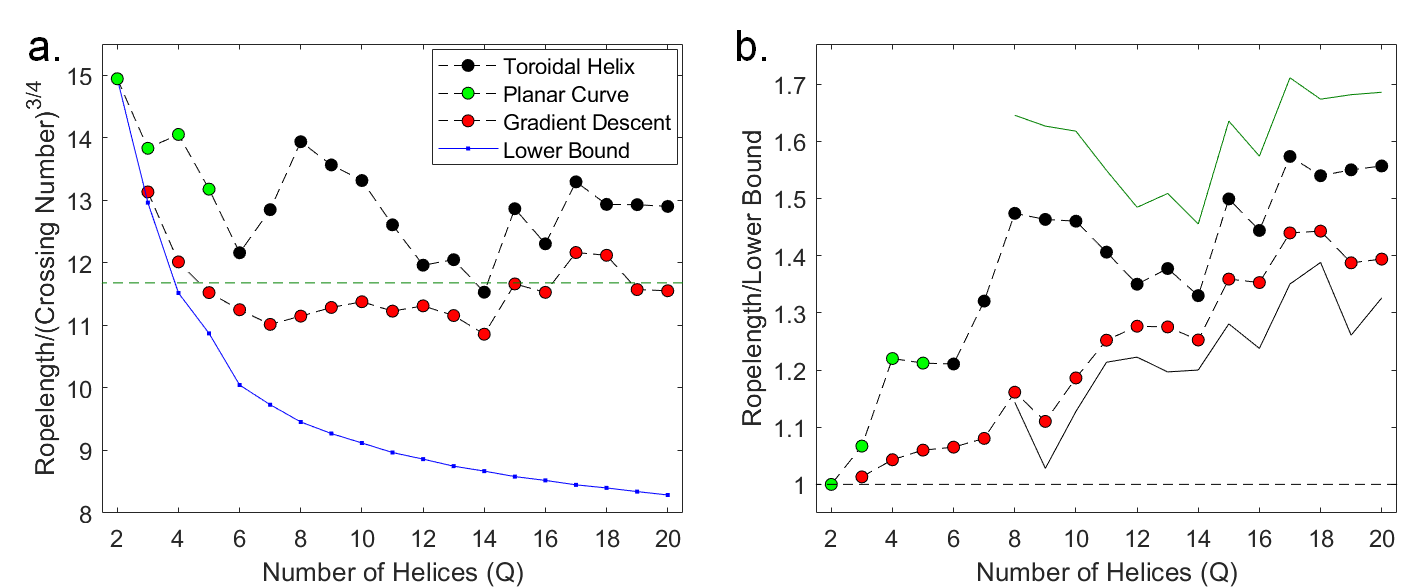}
    \caption{a. Three-quarter coefficients for the best constructed $T(Q,Q)$ links for $Q$ from 2 to 20, color-coded by their construction method. Values reached by gradient descent are shown in red and the lower bound is shown for comparison. The dashed line corresponds to the best large-$Q$ upper bound of 11.68. The red datum at $Q=3$ is from Ashton et al. \cite{ashton2011knot}. b. The same data as a ratio relative to the lower bound. The solid black line below the data indicates the ratio for the gradient descent data assuming that numerically optimized convex hulls minimal, and the green line indicates the ratio for the constructed links assuming the isoperimetric inequality. The large-$Q$ ratio of 1.77 sets the upper limit of this plot.}
    \label{fig:data}
\end{figure}

We will first discuss our results before describing how we arrived at them. We attempted to construct tight $T(Q,Q)$ links for $Q$ from 3 to 20 ($Q=2$ is trivial), and then used \textit{Ridgerunner}\cite{ashton2011knot} to further tighten each link. The data for $\alpha$ and the ratio to the lower bound are shown in Fig. \ref{fig:data}. Generally speaking, links could be constructed quite close to the lower bounds for $Q=3-6$. Beyond $Q=6$, the lower bound does not appear to be particularly constraining, with constructed links having an $\alpha$ between 12 and 13 (which continues out to $Q=60$), and tightened links with $\alpha$ between 11 and 12. The ratio with respect to the lower bound generally increases towards roughly 1.5 as the lower bound decreases. This leads us to suspect that the true bound for this type of link is closer to 10, and has much weaker $Q$-dependence. If the convex hulls found by Kallrath et al. \cite{kallrath2021near} are truly optimal, our links would move towards the lower bound by about 5 percentage points.


In addition to an interest in finding strong upper bounds, we have a desire to be able to simply describe optimal knots or links so that they may be easily reconstructed. For example, we recently demonstrated that the Mobius energy of a Hopf link is minimized when two congruent circles with perpendicular normals are displaced by $\sqrt{2}$ times their radius, a handy addendum to an earlier proof that the minimum is $2\pi^2$ \cite{agol}. A parametric equation with a small number of numeric parameters we deem more insightful than a list of one thousand Cartesian coordinates, although the latter is likely more efficient. In most cases it is not feasible to arrive at these parameters purely from geometric construction. Describing tight links with a small number of free parameters determined through numerical fitting may provide insight on the origins of these parameters and how they depend on crossing number, such that tight knots can be understood in greater depth in the future. In the next subsections we will describe our attempts to do this, and compare them to the results of numerical optimization. 


\subsection{Planar Curves}
Toroidal helices are not the only geometry from which torus links may be constructed, and the original goal of this paper was to find an alternative framework that would let us parameterize tighter links than those formed from concentric helices. This was ultimately unsuccessful, but we will described some insights gained along the way. The Mobius energy of a link is generally minimized by planar convex curves \cite{mobius}, and a $T(Q,Q)$ link can be constructed from congruent linked circles. These circles must be displaced and inclined relative to each other in order to link. We can construct the system as follows: start with $Q$ unit circles in the XZ plane centered on the origin. Rotate the $i$th circle about the Z axis by an angle $\phi=2\pi i/Q$. Translate each circle along the vector $\hat{r}_i=\langle\cos\phi,\sin\phi,0\rangle$ by a distance $\rho$, $0<\rho<1$. Rotate each circle, keeping its center fixed, about $\hat{r}_i$ by an angle $\psi$, $0<\psi<\pi/2$. The ropelength of this construction is $2\pi Q\frac{2}{d_{i,i+1}}$ where $d_{i,i+1}$ is the minimum distance between two sequential circles. The distance does not have a nice symbolic form but is easy to determine numerically. The ropelength can be minimized with respect to $\rho$ and $\psi$, and quickly converges to $\rho\approx1/2$ and $\psi\approx5\pi/18$ with ropelength $L\approx6.6C$, comparable to that of a $Q$-helix. Because this linear scaling is inefficient it was not deemed worthwhile to work out an exact expression. Equivalently, the circles may be arranged evenly in the XY plane with their centers on a circle of radius $\rho$ and tilted with angle of $\pi/2-\psi\approx2\pi/9$.

\begin{figure}
    \centering
    \includegraphics[width=1\linewidth]{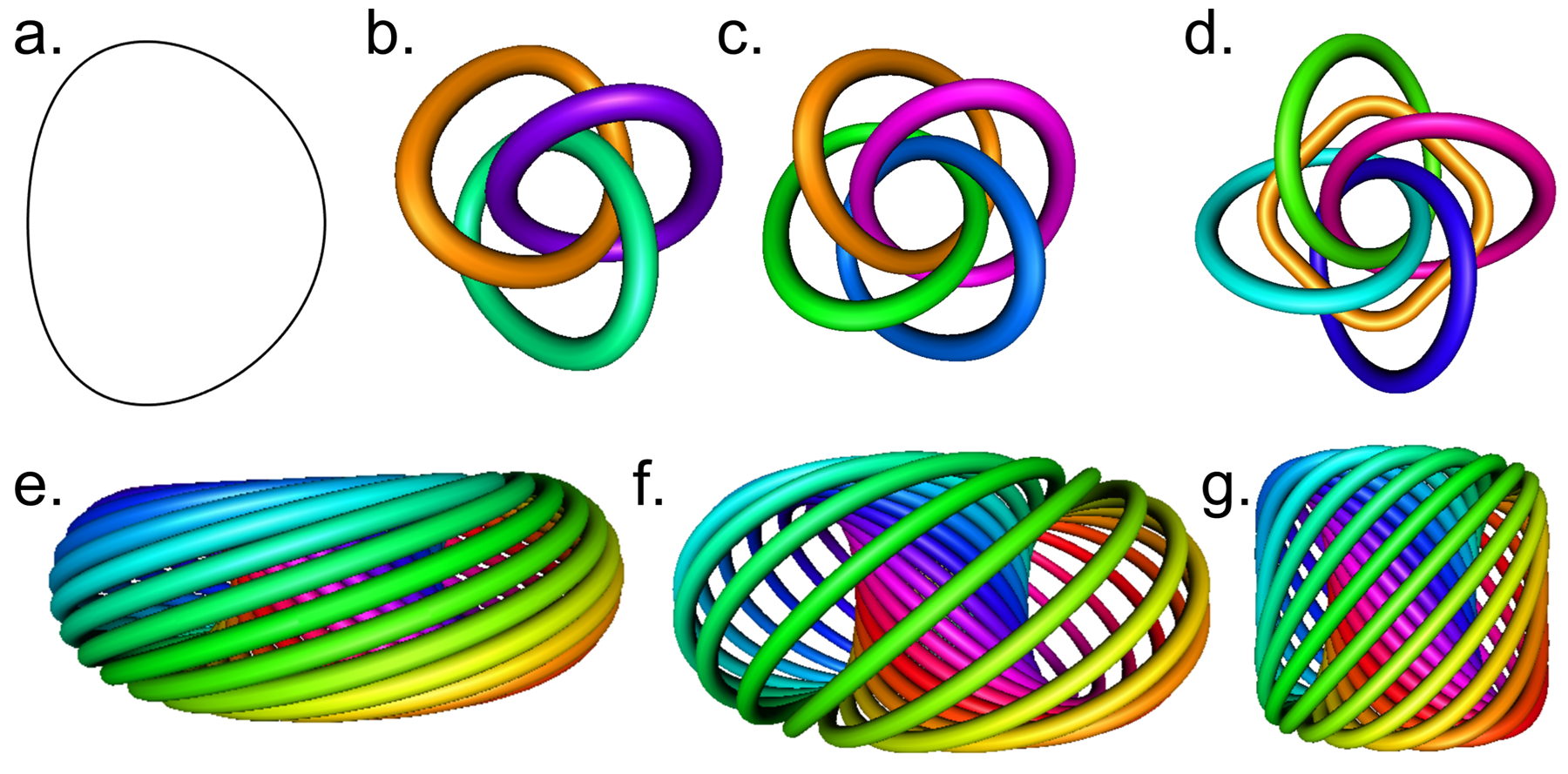}
    \caption{a. Typical shape of a gibbous curve that can be used to construct links. b-d. Optimal links for $Q=3,4,5$ constructed from gibbous curves, shown with reduced tube radius. d. includes a rounded square as the central curve. e-g. $T(20,20)$ links formed by optizing a 20-helix (e), congruent circles (f), and gibbous curves (g), shown with unit tube radius.}
    \label{fig:gibbous}
\end{figure}

Figure \ref{fig:gibbous}f shows a $T(20,20)$ constructed in such a manner. It is slightly narrower and taller than the equivalent 20-helix (Fig. \ref{fig:gibbous}e), but it is clear that there is excess length in the outer extremities of the circles. Deforming the circles into ellipses and using the eccentricity as a fit parameter will improve this slightly, but we found that curves resembling a gibbous moon (Fig. \ref{fig:gibbous}a) were even more efficient, as the portion of each curve that is far from the linked region can be made comparatively short. The parametric equation of a gibbous curve is:
\begin{equation}
    x(\theta)=\gamma\left(\cos(\theta)+\delta\cos(2\theta)\right),\ \ \ \  y=\sin(\theta),
\end{equation}
where the x-scale factor $\gamma$ is of order 1 and the Fourier coefficient $\delta$ is of order 0.1. This shape is convex for $\delta<1/4$. We can construct links similarly to our arrangement of circles, and minimize with respect to four variables: the displacement of each curve, the angle of inclination, the x-scale factor of the curves, and the Fourier coefficient. Examining a $T(20,20)$ link constructed from gibbous curves (Fig. \ref{fig:gibbous}g) shows an essentially cylindrical convex hull with a significantly smaller radius than the equivalent circular construction or 20-helix. Both $Q$-helices and linked circles have a ropelength that is asymptotically linear at approximately $L\approx6.6C$ The ropelength of these gibbous links is empirically $L\approx5.53C$. Since these scale linearly with crossing number they are inefficient for large links, but can be useful for small systems. We have constructed tight links by minimizing $\rho, \psi, \gamma$, and $\delta$ for $Q=3$ and $4$ (Fig. \ref{fig:gibbous}bc) and verifying that each curve had a minimum radius of curvature of at least 1. A tight $T(5,5)$ link was constructed from four gibbous curves linked through a central rounded square, where an additional two parameters, corresponding to the scale factor of the rounded square and the relative width of its straight segments, were used for fitting. It is possible to construct an efficient $T(5,5)$ with four circles and a rounded square as well, with only two free parameters. These are the three examples in which planar curves were found to be more efficient than toroidal helices, although it may be possible to extend this to higher $Q$. In the best case, the $T(3,3)$ constructed from gibbous curves was within 7\% of the lower bound, the others within 25\%.

\subsection{Toroidal Helices}
Our upper bound constructions may be inefficient for smaller systems because the major radius must be increased by a significant fraction in order to allow double linking. We constructed efficient toroidal helices by using the optimal shell arrangements found previously \cite{thompson2025ropelength}, initializing them with minor radii 2, 4, 6... but with the major radius as a free parameter for optimization. Compared to our previous methods of minimizing the total length subject to a no-overlap constraint, here we normalize the configurations by the minimum distance between any two helices, more akin to knot energy minimization. This allows, for example, 6 helices to be in the innermost shell for comparatively smaller systems. We then made small perturbations to the arrangement of the helices, including moving some to a perpendicular sub-link. 

Much like any linear $Q$-helix may be improved by smoothing one helix into a rod, a toroidal $T(Q,Q)$ link may be improved by taking an outer helix and transforming it into a circle with a normal perpendicular to the rest of the torus. More than one helix may be shifted to the perpendicular sub-link, and generally there is a benefit to doing so if it will reduce the number of shells in the main link. When a good configuration for a given $Q$ is found, an ideal configuration for $Q+1$ can be estimated by adding either another helix to the outer layer or forming a new layer, determining which has a shorter overall length, and then optimizing the arrangements by reverse Jenga. Although this does not necessarily yield the best configuration, it allows us to construct many efficient helices relatively quickly (typically a few seconds per configuration and tens of seconds to find an optimal configuration for each $Q$). We did this beyond $Q=60$ for $p=1$ and put some additional time into checking whether the arrangements were optimal for smaller $Q$, leading to the black points in Fig. \ref{fig:data}a. One particularly efficient link was $T(14,14)$, which was formed from two copies of a link formed by a circle surrounded by six helices, achieving $\alpha_{T_{14,14}}=11.53$. We carried out a similar procedure for $p=2$ and $3$, finding $\alpha$ typically between 11 and 12 up to $Q=60$, slightly lower than for $T(Q,Q)$ knots as expected. In the case of $p=3$ we typically found $\alpha$ below those presented previously \cite{thompson2025ropelength}, which were generated assuming the best helical arrangement is also the best toroidal arrangement.


\subsection{Gradient Descent}
Numerical knot tightening algorithms like \textit{Ridgerunner} \cite{ashton2011knot} and \textit{SONO} \cite{pieranski1998search} are the best tools for determining the upper bound of the ropelength of a given knot. The tightest $T(3,3)$ link ($6_3^3$/L6n1) found by Ashton et al. \cite{ashton2011knot} is only 1.3\% above the lower bound. Limitations of these methods include a high computational cost, sensitivity to initial conditions, uncertainty regarding convergence, and a difficulty in generalizing results beyond the study sample. Although it is not the primary goal of this work, we can use the knot tightening software \textit{Ridgerunner} to gain insight on what the tightest links look like and estimate our room for improvement, as well as to examine the relationship between the stronger upper bounds these algorithms produce and the lower bounds. We ran the first computations on our best constructed configurations, which are not guaranteed to be ideal initial conditions, with 200 vertices per component. Alternate computations were run using more symmetric initial conditions, which in some cases yielded superior results. $T(9,9)$ was an outlier with an  excess length until the vertex count was increased to 300 per component. At $Q=6$ and above, some ideal configurations adopted $Q-2$-fold symmetry, which we interpret as being due to the innermost and one of the outermost helices becoming minimized while the others conform to each other.

As seen in Fig. \ref{fig:data}, the \textit{Ridgerunner} ropelengths were about 10\% smaller than their constructed counterparts, and many are within 10\% of the lower bound. The results are typically below $\alpha_{opt}$ except at a few larger $Q$, which could likely be brought below with more computational attention. The ropelengths relative to the lower bounds typically increase towards about 1.4 with increasing $Q$, for which we believe the decreasing bound has a bigger effect than increasing inefficiency. These results do not necessarily imply that these values are close to the true ropelength, as the computations are sensitive to both initial conditions and local minima. For the purposes of documenting an improved upper bound for the universal $\alpha_0$, we report that we achieved low values for a few links that were beyond the scope of this work, the lowest being 10.19 achieved with a $T(20,10)$ link. Previously, we have achieved 10.02 with a $T(24,19)$ knot \cite{twitter}. Several failed attempts were made to find a knot or link with $\alpha<10$. We also have preliminary evidence that links constructed from a linking of a torus link and its mirror image cannot be tightened as much, with an estimated 15\% greater ropelength than the un-mirrored version.

\section{Conclusion}
Our goal was to construct knots for which the upper and lower ropelength bounds are relatively close. We have demonstrated this by constructing knots that are within a factor of 1.5 of the lower bound below 180 crossings, can be annealed to within 10\% of the lower bound for up to 42 crossings, and can be made within a factor of 1.77 for asymptotically large crossing number (requiring an assumption that the approximation in Eq. 20 holds). Improvements may be achieved by constructing toroidal helices with more complex harmonic modes, analogous to the non-circular curves used to construct efficient smaller links. Constructing links with higher modes suggests the idea of minimizing ropelength in the Fourier rather than Cartesian domain, which may prove a fruitful strategy for producing tight embeddings of knots. We achieved success in constructing tight links by creating a Hopf link of two smaller links. Such a construction could be performed with other knots or links, e.g. winding a dense helical cylinder into a non-alternating torus knot. This may prove a fruitful avenue for further investigation of ropelength bounds.

\section{Acknowledgments}
ARK is supported by the National Science Foundation, grant number 2336744. Ridegerunner computations were performed on the University of Edinburgh Research Compute Cluster. The author is grateful to Kathryn McCormick, Kasturi Barkataki, and Finn Thompson for insightful discussions. After the first version of this paper was uploaded to arXiv, the author was contacted by Terence Tao with a link to ChatGPT's analysis of the preprint. Some of the feedback from ChatGPT was incorporated into the following version. This includes more clearly stated assumptions underlying some main results, correcting several transcription errors, and updating the equation for and table of helical corrections in the Appendix.

\bibliographystyle{unsrt}
\bibliography{knotrefs}

\section{Appendix}
\subsection{Toroidal Correction}
Table 1, which is hopefully nearby, shows the ratio of the arc length of a toroidal helix with a given ratio of major to minor radius, to that of a straight helix with a height such that $2\pi R=pH$. Specifically it is computing:
\begin{equation} 
\frac{\int_0^{2\pi}\sqrt{r^2+p^{-2}(R-r\cos\theta)^2}d\theta}{\sqrt{(2\pi R/p)^2+(2\pi r)^2}}=\frac{\int_0^{2\pi}\sqrt{1+p^{-2}(R/r-\cos\theta)^2}d\theta}{2\pi \sqrt{( R/rp)^2+1}}
\end{equation}
as a function of $R/r$. This was computed by generating toroidal helices with 10,000 vertices and measuring the arc length discretely. Data are shown for $p=1,2,3$. The excess decreases with the negative fourth power of the ratio. The correction is smaller for larger $p$ when the ratio is small, but larger for larger $p$ when the ratio is large. Note that with a ratio of 1 a toroidal helix is self-intersecting for $p>1$.

\begin{table}[]
\begin{tabular}{|l|l|l|l|l|l|l|l|}
\hline
Ratio & $p=1$      & $p=2$      & $p=3$      & Ratio & $p=1$       & $p=2$      & $p=3$      \\ \hline
1     & 1.06998  & 1.039895 & 1.022286 & 3     & 1.0026752 & 1.006131 & 1.007052 \\ \hline
1.1   & 1.058092 & 1.036947 & 1.021417 & 3.5   & 1.001501  & 1.003911 & 1.005065 \\ \hline
1.2   & 1.04807  & 1.034065 & 1.020521 & 4     & 1.000902  & 1.002572 & 1.00366  \\ \hline
1.3   & 1.039731 & 1.031288 & 1.019608 & 4.5   & 1.0005727 & 1.001743 & 1.002672 \\ \hline
1.4   & 1.032854 & 1.028642 & 1.018687 & 5     & 1.0003803 & 1.001215 & 1.001975 \\ \hline
1.5   & 1.027218 & 1.026147 & 1.017768 & 5.5   & 1.000262  & 1.000869 & 1.001481 \\ \hline
1.6   & 1.022615 & 1.023815 & 1.016856 & 6     & 1.0001863 & 1.000635 & 1.001125 \\ \hline
1.7   & 1.018859 & 1.021651 & 1.015959 & 6.5   & 1.0001359 & 1.000474 & 1.000866 \\ \hline
1.8   & 1.015793 & 1.019657 & 1.015083 & 7     & 1.0001015 & 1.000361 & 1.000676 \\ \hline
1.9   & 1.013286 & 1.017828 & 1.014232 & 7.5   & 1.0000773 & 1.000279 & 1.000534 \\ \hline
2     & 1.011229 & 1.016158 & 1.013409 & 8     & 1.0000598 & 1.000218 & 1.000426 \\ \hline
2.1   & 1.009536 & 1.014639 & 1.012617 & 8.5   & 1.0000471 & 1.000174 & 1.000344 \\ \hline
2.2   & 1.008137 & 1.013262 & 1.011859 & 9     & 1.0000375 & 1.00014  & 1.00028  \\ \hline
2.3   & 1.006976 & 1.012016 & 1.011136 & 9.5   & 1.0000303 & 1.000113 & 1.00023  \\ \hline
2.4   & 1.006008 & 1.010891 & 1.010448 & 10    & 1.0000247 & 1.000093 & 1.000191 \\ \hline
2.5   & 1.005197 & 1.009876 & 1.009795 & 20    & 1.0000015 & 1.000006 & 1.000013 \\ \hline
 
\end{tabular}
\caption{Table of toroidal corrections for the length of a helix as a function of the major to minor radius ratio.}
\end{table}

\subsection{Validity of the Rectangular Approximation}
Here we present a sketch of a proof that the non-overlapping helices can be constructed by using one less helix than proscribed by the small-angle approximation. This involves some unwieldy expressions that were determined using \textit{Mathematica} and may not be optimally simplified.
The squared distance between two adjacent helices on a cylinder is:
\begin{equation}
    d^2(\theta_\Delta,N,r,h)=2r^2\left(1-\cos\left(\theta_\Delta+\frac{2\pi}{N}\right)\right)+h^2\theta_\Delta^2
\end{equation}
This is minimized by $\theta_{\Delta o}$, and we require that $d^2(\theta_{\Delta o})\geq4$, but this equation is transcendental and we cannot derive a symbolic solution. We desire the maximum value of $N$ that allows a solution for a given $r$ and $h$. When $N$ is large the angle is small, and if we either approximate the cosine as a quadratic or find the closest points between two lines parallel to the diagonal of an equivalent cylinder, we find an approximate maximum $N$:
\begin{equation}
    N_{a}=\frac{\pi hr}{\sqrt{h^2+r^2}}
\end{equation}
We need to demonstrate that $d^2(N_a)\geq4$ in order to show that that many helices can fit around a cylinder. Here we sketch an argument that the constraint is satisfied when $N_a-\epsilon$ cylinders are used. First, we expand $d^2(N_a-\epsilon)$ as a Taylor series in $\theta_\Delta$ and truncate it at the cubic term:
\small
\begin{equation}
    d_c^2(\theta_\Delta,h,r)=-\frac{1}{3} r^2 \theta_\Delta ^3 \sin \left(\frac{2 \pi }{\frac{\pi h r}{\sqrt{h^2+rr^2}}-\epsilon}\right)+\theta_\Delta ^2 \left(r^2 \cos \left(\frac{2 \pi }{\frac{\pi h r}{\sqrt{h^2+r^2}}-\epsilon}\right)+h^2\right)+2 r^2 \theta_\Delta \sin \left(\frac{2 \pi }{\frac{\pi h r}{\sqrt{h^2+r^2}}-\epsilon}\right)+4 r^2 \sin ^2\left(\frac{\pi }{\frac{\pi h r}{\sqrt{h^2+r^2}}-\epsilon}\right)
\end{equation}
\normalsize
The derivative with respect to $\theta_\Delta$ has two zeros, the closest of which to $\theta_\Delta=0$ is relevant:
\small
\begin{equation}
\theta_{\Delta co}= -\frac{\csc \left(\frac{2 \pi }{\frac{\pi h r}{\sqrt{h^2+r^2}}-\epsilon}\right) \left(\sqrt{4 \left(r^2 \cos \left(\frac{2 \pi }{\frac{\pi h r}{\sqrt{h^2+r^2}}-\epsilon}\right)+h^2\right)^2+8 r^4 \sin ^2\left(\frac{2 \pi }{\frac{\pi h r}{\sqrt{h^2+r^2}}-\epsilon}\right)}-2 \left(r^2 \cos \left(\frac{2 \pi }{\frac{\pi h r}{\sqrt{h^2+r^2}}-\epsilon}\right)+h^2\right)\right)}{2 r^2}
\end{equation}
\normalsize
Substituting this into $d_c^2$ we find the minimum approach of the cubic approximation:

\small
\begin{multline}
     d^2_{co}(r,h)= 4 r^2 \sin ^2\left(\frac{\pi }{\frac{\pi h r}{\sqrt{h^2+rr^2}}-\epsilon}\right)-\zeta+\frac{\csc ^2\left(\frac{2 \pi }{\frac{\pi h r}{\sqrt{h^2+rr^2}}-\epsilon}\right) \zeta^3}{24 r^4}+\frac{\left(r^2 \cos \left(\frac{2 \pi }{\frac{\pi h r}{\sqrt{h^2+rr^2}}-\epsilon}\right)+h^2\right) \csc ^2\left(\frac{2 \pi }{\frac{\pi h r}{\sqrt{h^2+rr^2}}-\epsilon}\right) \zeta^2}{4 r^4},\\ \zeta=\sqrt{4 \left(r^2 \cos \left(\frac{2 \pi }{\frac{\pi h r}{\sqrt{h^2+rr^2}}-\epsilon}\right)+h^2\right)^2+8 r^4 \sin ^2\left(\frac{2 \pi }{\frac{\pi h r}{\sqrt{h^2+rr^2}}-\epsilon}\right)}-2 \left(r^2 \cos \left(\frac{2 \pi }{\frac{\pi h r}{\sqrt{h^2+rr^2}}-\epsilon}\right)+h^2\right)     
\end{multline}
\normalsize

We need to demonstrate that this is always greater than 4 for relevant values of of $r$ and $h$, e.g. that $r>2$, $r<h$. This is done by examining several limits, in the tall cylinder limit as $h$ goes to infinity, and in the isotropic cases when $r=h$ for $r=2$ and $h\rightarrow\infty$. Starting with the simple example of using one less than the approximate number of helices, $\epsilon=1$, in the tall cylinder case we find:
\begin{equation}
    \lim_{h\rightarrow\infty}d^2_{0c}(h,r)=4r^2\sin^2\frac{\pi}{\pi r-1}\approx4+\frac{8}{\pi r}
\end{equation}
This is always greater than 4, tending to 4 at $r\rightarrow\infty$ and 5.2 when $r=2$. In the isotropic case we have set $r=h$ and have:
\small

\begin{multline}
   d^2_{0c}(h)= 4 h^2 \sin ^2\left(\frac{2 \pi }{\sqrt{2} \pi \sqrt{h^2}-2}\right)+4 h^2 \cos ^2\left(\frac{2 \pi }{\sqrt{2} \pi \sqrt{h^2}-2}\right)-2 \sqrt{2} \sqrt{-h^4 \cos ^2\left(\frac{2 \pi }{\sqrt{2} \pi \sqrt{h^2}-2}\right) \left(\cos \left(\frac{4 \pi }{\sqrt{2} \pi \sqrt{h^2}-2}\right)-3\right)}\\+\left(\frac{1}{2h^2}-\frac{1}{3h^4}\right)\left(h^2+h^2 \cos \left(\frac{4 \pi }{\sqrt{2} \pi \sqrt{h^2}-2}\right)-\sqrt{2} \sqrt{-h^4 \cos ^2\left(\frac{2 \pi }{\sqrt{2} \pi \sqrt{h^2}-2}\right) \left(\cos \left(\frac{4 \pi }{\sqrt{2} \pi \sqrt{h^2}-2}\right)-3\right)}\right)^3 \csc ^2\left(\frac{4 \pi }{\sqrt{2} \pi \sqrt{h^2}-2}\right)
\end{multline}
\normalsize
When $h=2$, this has an evaluation that approximates to 6.38. In the other limit the distance tends to 4 and we have:
\begin{equation}
        \lim_{h\rightarrow\infty}d^2_{0c}(R_0)\approx4+\frac{8\sqrt{2}}{\pi h}
\end{equation}
If instead of subtracting 1 we wish to find the smallest allowable value of $\epsilon$, we can examine the most constraining case of a narrow cylinder and large height, and find that:

\begin{equation}
    \pi\left(r-\left(\arcsin\frac{1}{r}\right)^{-1}\right)\leq\epsilon\leq2\pi-6
\end{equation}

The minimal value is reached at $r=2$ when $N_a=2\pi$ and the actual maximum is 6. It is desirable to show that the constraint is satisfied for $\epsilon=0$ when $N=N_{a}$, however, following the same arguments we find in the tall cylinder limit:
\begin{equation}
    d_c^2(r)\approx4-\frac{4}{3r^2}<4
\end{equation}
and in the isotropic limit
\begin{equation}
    d_c^2(h)\approx4-\frac{68}{45h^4}<4
\end{equation}
This suggests that the constraint is not satisfied for $\epsilon=0$. This is the extent of our proof, which is incomplete because we have not shown that the full constraint equation is satisfied, only the cubic approximation. In all relevant limits, the cubic approximation of the constraint equation is satisfied by one less than the maximum number of helices gleaned from the rectangular approximation. Since the number of helices is an integer, using $N=\lfloor N_a\rfloor$ in practice would likely not cause constraint violations with frequency, but we cannot prove that it will work in all cases. We can, however, show numerically that when $N_a-(2\pi-6)$ and $r=2$ are substituted into the full constraint equation, the constraint is satisfied for all tested $h$, with the limiting squared distance being $\approx4+6.2/h^2$.

\end{document}